\documentclass[11pt]{amsart}
\usepackage{latexsym,amssymb,amsmath,amscd,amsthm}
\numberwithin{equation}{section}
\theoremstyle{remark}
\newtheorem{theorem}{{\bf THEOREM}}[section]

\newcommand{\bq}{\begin{equation}}
\newcommand{\bea}{\begin{array}}
\newcommand{\eea}{\end{array}}

\newcommand{\mf}{\mathfrak}

\newcommand{\gG}{\Gamma}
\newcommand{\gt}{\theta}
\newcommand{\gs}{\sigma}

\newcommand{\gag}{\gamma}
\newcommand{\gd}{\delta}
\newcommand{\pp}{\partial}

\newcommand{\tl}{\tilde}

\newcommand{\bs}{\blacksquare}

\newcommand{{\DDD}}{D\!\!\!\!\!\!-}

\title{TOWARD A CANONICAL QKDV EQUATION}
\author{Robert Carroll\\University of Illinois, Urbana, IL 61801}

\date{March, 2003\thanks{email: rcarroll@math.uiuc.edu}}

\begin{document}

\bibliographystyle{plain}

\maketitle


\section{INTRODUCTION}
\renewcommand{\theequation}{1.\arabic{equation}}
\setcounter{equation}{0}

In \cite{c1,c2} we gave a number of formulas for qKdV equations derived in a
q-Virasoro context.  The constructions were meaningful even though the cocycle terms
did not visibly satisfy a Jacobi identity.  In the present note we modify the 
framework slightly and produce a cocycle term (morally equivalent to one of those
used in \cite{c1,c2}) and thereby exhibit a more convincing candidate for a canonical
qKdV equation. It remains to see how this is connected to the qKdV equation from the
hierarchy picture. For background on q-analysis etc. we refer to \cite{c3}.

\section{BACKGROUND}
\renewcommand{\theequation}{2.\arabic{equation}}
\setcounter{equation}{0}

We recall first from \cite{c1,c2,l1,l2} the following information regarding q-Virasoro
constructions.
Thus
work on $S^1$ with ($q\ne 0,\pm 1$)
\bq\label{54}
\pp_qz=\frac{q^mz^m-q^{-m}z^m}{(q-q^{-1}z}=z^{m-1}[m];\,\,[m]=
\frac{q^m-q^{-m}}{q-q^{-1}}
\end{equation}
We adapt the formalism of \cite{l1} as follows.  Let $D_n=-z^{n+1}\pp$ with
$\pp:\,z^m\to q^m[m]z^{m-1}$ so $\pp\sim\pp_q\tau$ where $\tau f(z)=f(qz)$.
Generally we will think of $z=e^{i\gt}\in S^1$ so $(1/2\pi
i)\int_{S^1}z^ndz=(1/2\pi)\int z^{n+1}d\gt=\gd_{(-1,0)}$ which will be written
as ${\bf (A1)}\,\,\int z^n=\gd_{(-1,0)}$.  Write also ${\bf (A2)}\,\,
\ell_n\sim D_n=-z^{n+1}\pp_q\tau$ (we used $D_n\sim-\ell_n$ in
\cite{c1,c2} - this produces a few sign changes there but does not affect any
basic conclusions).   It is
known that q-brackets are needed now where
\bq\label{55}
[\ell_m,\ell_n]_q=q^{m-n}\ell_m\ell_n-q^{n-m}\ell_n\ell_m=[m-n]\ell_{m+n}
\end{equation}
For a central term in a putative $Vir_q$ one wants (cf. \cite{c1,c2,l1,l2})
a formula ${\bf (A3)}\,\,c[m+1][m][m-1]\gd_{m+m,0}$ (see below for an optimal
term). First we want to
formulate the q-bracket in terms of vector fields as follows (the central
term will be added later in a somewhat ad hoc manner).  This can be done as a
direct calculation using the basic definition of $\pp$ above (cf. also \eqref{C} below).
Thus
\bq\label{56}
[z^n\pp,z^m\pp]_q\sim q^{n-m}z^n\pp(z^m\pp)-q^{m-n}z^m\pp(z^n\pp)=
\end{equation}
$$=(q^n[m]-q^m[n])z^{m+n-1}\pp=
[m-n]z^{m+n-1}\pp=[n-m](-z^{m+n-1}\pp)$$
Let now $v\sim \sum a_nz^n$ and $w\sim \sum b_mz^m$; then we define a bracket in
$Vec(S^1)$ via
\bq\label{57}
[v\pp,w\pp]_q=-\sum a_nb_m
[n-m]z^{m+n-1}\pp
\end{equation}
We defined a bracket of vector fields in \cite{c1,c2} so that from \eqref{57} there
resulted a correspondence
\bq\label{58}
v'w-vw'\sim -[v\pp_x,w\pp_x]\sim -[v\pp,w\pp]_q=-\{(\tau v)(\pp_qw)-(\tau
w)(\pp_qv)\}\tau
\end{equation}
This dangling $\tau$ created some complications in further calculation however and is
removed below in the new formulation.
\\[3mm]
\indent
{\bf REMARK 2.1.}
In \cite{l1} one defines the q-analogue of the enveloping algebra 
of the Witt algebra ${\mf W}$ as the associative algebra ${\mf U}_q({\mf W})$
having generators $\ell_m\,\,(m\in{\bf Z})$ and relations \eqref{55}.  The q-deformed
Virasoro algebra is defined as the associative algebra ${\mf U}_q(Vir)$ having generators
$\ell_m\,\,(m\in {\bf Z})$ and relations ($q\ne$ root of unity)
\bq\label{59}
q^{m-n}\ell_m\ell_n-q^{n-m}\ell_n\ell_m=[m-n]\ell_{m+n}+\gd_{m+n,0}\frac
{[m+1][m][m-1]}{[2][3]<m>}\hat{c}
\end{equation}
where $<m>=q^m+q^{-m}$ and $\hat{c}\ell_m=q^{2m}\ell_m\hat{c}$ (thus $\hat{c}$ is an
operator which we examine below and we refer to \cite{c1,c2,l1,l2} for the central
term).   Then
${\mf U}_q(Vir)\sim Vir_q$ is a {\bf Z} graded algebra with $deg(\ell_m)=m$ and
$deg(\hat{c})=0$.  One also introduces in 
\cite{l1} a larger algebra ${\mf U}(V_q)=$ associative algebra generated by $J^{\pm 1},
\,\hat{c},\,d_m\,\,(m\in{\bf Z})$ with relations
\bq\label{60}
JJ^{-1}=J^{-1}J=1;\,\,Jd_mJ^{-1}=q^md_m;\,\,\hat{c}J=J\hat{c};\,\,\hat{c}d_m=q^md_m\hat{c};
\end{equation}
$$q^md_md_nJ-q^nd_nd_mJ=[m-n]d_{m+n}+\gd_{m+n,0}\frac{[m+1][m][m-1]}{[2][3]<m>}\hat{c}$$
The subalgebra of ${\mf U}(V_q)$ generated by $\ell_m'=d_mJ$ and
$\hat{c}'=\hat{c}J\,\,(m\in {\bf Z})$ is the same as ${\mf U}_q(Vir)$. It is stated in 
\cite{l2} that $V_q$ is the universal quantum central extension of ${\mf W}_q$
and thus \eqref{60} is better
adapted for optimal algebraic and geometric meaning; it is this aspect which we
emphasize in this paper.
$\hfill\bs$
\\[3mm]\indent
In \cite{c1,c2} we constructed various forms of qKdV based on $\ell_m$ and 
\eqref{59} but will now redo
some constructions in terms of $d_m,\,\hat{c},$ and $J^{\pm 1}$ in hopes of producing
a genuine cocycle.  Recall we had constructed various pseudo-cocycles in the form
(omitting multiplicative factors of $q$)
\bq\label{A}
\psi(f\pp,g\pp)=\int(\tau\pp^3f)(\tau g);\,\,\psi'=\int (\tau g)\tau
(\pp^2\tau^2(\tau+\tau^{-1})^{-1}\pp\tau^{-5}f;
\end{equation}
$$\psi''(f\pp,g\pp)=\int(\tau g)(\tau^{-4}\pp^2(\tau+\tau^{-1})^{-1}\pp(\tau^2f);\,\,
\psi'''(f\pp,g\pp)=\int (\tau g)(\tau(\tau^{-5}\pp^3\tl{\gG}^{-1}f)$$
where $\tl{\gG}=(q^{-1}\tau+q\tau^{-1})$.  We recall also from \cite{l2}
\bq\label{B}
[d_m,d_n]=Jd_mJ^{-1}d_nJ-Jd_nJ^{-1}d_mJ=
\end{equation}
$$=q^md_md_nJ-q^nd_nd_mJ=q^{-n}Jd_m
d_n-q^{-m}Jd_nd_m$$
Further $\hat{c}d_m=q^md_m\hat{c}$ suggests $\hat{c}=\tau$ here.  Indeed we first 
correct the definition of $d_n$ from \cite{l2} since $d_m$ is being used as $\ell_m
\tau^{-1}$.  Dropping the minus sign momentarily, from $\ell_m=z^{m+1}\pp_q\tau$ we get
then
${\bf (A4)}\,\,d_m=z^{m+1}\pp_q$. Then, using ${\bf (A5)}\,\,\pp_q\tau=q\tau\pp_q$, one
obtains $\tau d_m=\tau(z^{m+1}\pp_q)=q^mz^{m+1}
\pp_q\tau$ and $\tau^{-1}d_m=q^{-m}d_m\tau^{-1}$. 
In addition $Jd_mJ^{-1}=q^md_m$ corresponds to $Jd_m=q^md_mJ$ so we identify
$J=\tau$.  Writing $\ell_m=d_mJ=d_m\tau$ we can also easily see that the brackets
$[d_m,d_n]$ above are exactly the q-brackets ${\bf (A6)}\,\,[\ell_m,\ell_n]_q=
q^m\ell_m\ell_n-q^n\ell_n\ell_m$.
\\[3mm]\indent
Now in \cite{l2} a Jacobi type identity is used involving an operator $\gs(x)=
(1/2)(\tau+\tau^{-1})(x)$ for $x\in \oplus {\bf C}d_n$.  This seems to be better phrased
in terms of an operator ${\bf (A7)}\,\,\gG(d_p)=<p>d_p$ which avoids the need to 
carry $\tau$ around otherwise.  Then we can check that the rule in \eqref{60}, rewritten
as ${\bf (A8)}\,\,[d_m,d_n]=[m-n]d_{m+n}+\gag_m\gd_{m+n,0}\hat{c}$ will yield
\bq\label{P}
[[d_m,d_n],\gG(d_p)]+[[d_n,d_p],\gG(d_m)]+[[d_p,d_m],\gG(d_n)]=\Xi_{m,n,p}=0
\end{equation}
This is based on two identities; one, stated in \cite{c1}, is
\bq\label{Z}
[m-n][m+n-p]<p>+[n-p][n+p-m]<m>+[p-m][p+m-n]<n>=0
\end{equation}
The second is
\bq\label{X}
[p+1][p][p-1][m-n]+[m+1][m][m-1][n-p]+
\end{equation}
$$+[n+1][n][n-1][p-m]=0;\,\,(m+n+p=0)$$
The first proof is straightforward and for the second we note that one can write
it in the form
\bq\label{Y}
[m+1][m][m-1][2n+m]-[m+n-1][m+n][m+n+1][m-n]-
\end{equation}
$$-[n+1][n][n-1][n+2m]=0$$
we recall $<p>=q^p+q^{-p}$ and the trick is to rewrite \eqref{Y} in the form
\bq\label{W}
[(<2m+1>-<1>)(<2n+2m-1>-<2n+1>)]-
\end{equation}
$$-[(<2m+2n-1>-<1>)(<2m+1>-<2n+1>)]-$$
$$-[(<2n+1>-<1>)(<2n+2m-1>-<2m+1>)]=0$$
(the 0 following by a simple calculation).  
Now to prove \eqref{P} we write e.g.
\bq\label{V}
[[d_m,d_n],<p>d_p]=[[m-n]d_{m+n}+\gag_m\gd_{m+n,0}\hat{c},<p>d_p]
\end{equation}
and note that $\hat{c}=\tau\sim d_0$ so from \eqref{B} $[\hat{c},d_p]=q^0\tau d_p\tau-
q^pd_p\tau^2=q^p(1-1)(d_p\tau^2)=0$.  Hence we get
\bq\label{U}
\Xi_{m,n,p}=[[m-n]d_{m+n},<p>d_p]+
\end{equation}
$$+[[n-p]d_{n+p},<m>d_m]+[[p-m]d_{p+m},<n>d_n]$$
Reversing the bracket order, there will be terms
\bq\label{T}
\left([m-n]<p>[p-m-n]+[n-p]<m>[m-n-p]+\right.
\end{equation}
$$\left.+[p-m]<n>[n-p-m]\right)d_{p+m+n}$$
which vanishes by \eqref{Z}, and 
\bq\label{S}
\left([m-n]<p>\gag_p+[n-p]<m>\gag_m+[p-m]<n>\gag_n\right)\gd_{p+m+n,0}\hat{c}
\end{equation}
which is zero by \eqref{X}.  This shows that \eqref{P} will hold and $V_q$ will be a
genuine central extension of ${\mf W}_q$, with a reasonable Jacobi identity \eqref{P}.

\section{COCYCLES}
\renewcommand{\theequation}{3.\arabic{equation}}
\setcounter{equation}{0}

The cocycle search involved finding $\psi(f\pp,g\pp)$ where $\pp\sim\pp_q\tau$ and we
gave some pseudo-examples $\psi',\,\psi'',\psi'''$ in \eqref{A}.  Now in \eqref{56}
- \eqref{58} we recall $\ell_m\sim -z^{m+1}\pp=-z^{m+1}\pp_q\tau$ and a $d_m$
formulation would drop the $\tau$.  Thus work with $\pp_q$ instead of $\pp=\pp_q\tau$
with ${\bf (A9)}\,\,\pp_qz^{p+1}=q^{p+1} z^{p+1}\pp_q+[p+1]z^p\tau^{-1}$ based on
$\pp_qf=(\tau f)\pp_q+(\pp_qf)\tau^{-1}$.  Then
\bq\label{C}
[d_m,d_n]=[z^{m+1}\pp_q,z^{n+1}\pp_q]=q^md_md_n\tau-q^nd_nd_m\tau=q^mz^{m+1}
\pp_qz^{n+1}\pp_q\tau-
\end{equation}
$$-q^nz^{n+1}\pp_qz^{m+1}\pp_q\tau=[n-m]z^{n+m+1}\pp_q=[m-n]d_{m+n}$$
In this context the $\gd_{m+n,0}$ term does not arise.  Note here $\tau^{-1}\pp_q=
q\pp_q\tau^{-1}$ and
${\bf (A10)}\,\,
q^{m+1}[n+1]-q^{n+1}[m+1]=[n-m]$.
Let now $v\sim \sum v_{n+1}z^{n+1}$ and 
$w=\sum w_{m+1}z^{m+1}$; then
\bq\label{E}
[v\pp_q,w\pp_q]=[\sum v_{n+1}d_n,\sum w_{m+1}d_m]=\sum v_{n+1}w_{m+1}[d_n,d_m]=
\end{equation}
$$=\sum v_{n+1}w_{m+1}[n-m]z^{n+m+1}\pp_q$$
Going back to \eqref{C} this corresponds then to
\bq\label{F}
[v\pp_q,w\pp_q]=[(\tau v)(\pp_qw)-(\tau w)(\pp_qv)]\pp_q
\end{equation}
and this is exactly \eqref{58} but with the offending $\tau$ removed.
\\[3mm]\indent
Now try to build in a cocycle term automatically by using \eqref{60} in \eqref{C}
so that a term arises of the form
\bq\label{G}
[z^{m+1}\pp_q,z^{n+1}\pp_q]=[d_m,d_n]=[m-n]d_{m+n}+\gag_m\gd_{m+n,0}\tau
\end{equation}
For $[v\pp_q,w\pp_q]$ we get then an additional term
\bq\label{H}
\sum v_{n+1}w_{m+1}\gag_n\gd_{m+n,0}\tau=\sum v_{n+1}w_{-n+1}
\frac{[n+1][n][n-1]\tau}{<n>[2][3]}
\end{equation}
Now consider integrals involving $w$ and
\bq\label{I}
(\pp_q^3v)=\sum v_{n+1}[n+1][n][n-1]z^{n-2}
\end{equation}
over $S^1$.  If we write $(\tau+\tau^{-1})z^n=(q^n+q^{-n})z^n=<n>z^n$ then look at
\bq\label{J}
(\pp_q^2(\tau+\tau^{-1})^{-1}(\pp_qv))=\sum v_{n+1}\frac{[n+1][n][n-1]}{<n>}z^{n-2}=
\sum v_{n+1}\tl{\gag}z^{n-2}
\end{equation}
This starts to resemble $\psi''$ in \eqref{B}.  Then for $a=1/[2][3]$ 
\bq\label{K}
a\int w(\pp_q^2(\tau+\tau^{-1})^{-1}(\pp_qv))=\int\sum w_{m+1}v_{n+1}z^{m+n-1}
\gag_n=
\end{equation}
$$=\sum v_{n+1}w_{-n+1}\gag_n=\phi(v\pp_q,w\pp_q)$$
which agrees with \eqref{A}.  Consider $\phi(w\pp_q,v\pp_q)$ with, for $n\to -n$,
${\bf (A11)}\,\,\sum w_{n+1}v_{-n+1}\gag_n\tau=-\sum v_{n+1}w_{-n+1}\gag_n\tau$ since 
$\gag_{-n}=-\gag_n$.  Hence $\phi$ is antisymmetric.  For a Jacobi condition we  
go to \eqref{P} and consider (cf. \cite{c1,c2})
\bq\label{L}
\phi([v\pp_q,w\pp_q],\gG(u\pp_q))+\phi([w\pp_q,u\pp_q],\gG(v\pp_q))+
\phi([u\pp_q,v\pp_q],\gG(w\pp_q))=\Upsilon(u,v,w)
\end{equation}
Since $\phi(v\pp_q,w\pp_q)=-\phi(w\pp_q,v\pp_q)$ we can reverse all brackets and
consider
\bq\label{M}
\gG(u\pp_q)\sim \sum u_{p+1}<p>z^{p+1}\pp_q;\,\,[v\pp_q,w\pp_q]\sim \sum
v_{n+1}w_{m+1}[n-m]z^{n+m+1}\pp_q;
\end{equation}
$$a\pp_q^2(\tau+\tau^{-1})^{-1}\pp_qu=\sum u_{p+1}\gag_pz^{p-2}$$
Then \eqref{K} gives
\bq\label{1}
\phi([v\pp_q,w\pp_q],\gG(u\pp_q))=-\phi(\gG(u\pp_q),[v\pp_q,w\pp_q])=
\end{equation}
$$=
-\int\sum<p>u_{p+1}\gag_pz^{p-2}\sum v_{n+1}w_{m+1}[n-m]z^{n+m+1}=$$
$$=-\sum<p>\gag_p[n-m]u_{p+1}v_{n+1}w_{m+1}\gd_{m+n+p,0}$$
Similarly for each $u_{p+1}v_{n+1}w_{m+1}$ we will have contributions with coefficients
$[n-p]<m>\gag_m$ and $[p-m]<n>\gag_n$ (cf. \eqref{S}), the sum of which is zero by
\eqref{X}.
Consequently
\begin{theorem}
The term $\phi(v\pp_a,w\pp_q)$ in \eqref{K} is a cocycle and following the constructions
in \cite{c1,c2} one has a possibly canonical qKdV equation in the form ($a=1/[2][3]$)
\bq\label{2}
u_t+c'\pp_q^2(\tau+\tau^{-1})^{-1}\pp_qu+\pp_q(u\tau u)+\tau^{-1}u\pp_q\tau^{-1}u
\end{equation}
\indent
{\it Proof.}  We modify slightly the constructions in \cite{c1,c2} and take a duality
expression
\bq\label{3}
<(v\pp_q,a),(u,c)>=\int vu+ac
\end{equation}
(cf. \cite{a1}).  Then write
\bq\label{4}
q<[f\pp_q,a),(g\pp_q,b)],(u,c)>=-q\int [f\pp_q,g\pp_q]u+cq\phi(f\pp_q,g\pp_q)=
\end{equation}
$$=-q\int [(\tau f)(\pp_qg)-(\tau f)(\pp_qf)]u+caq\int g\pp_q^2(\tau+\tau^{-1})^{-1}
\pp_q^2f$$
We note from \cite{c1} that ${\bf (A12)}\,\,q\int fg=\int \tau^{-1}f\tau^{-1}g$ and
$\int\pp_qf=0$ while from ${\bf (A13)}\,\,\pp_q(gfu)=\pp_qf\tau^{-1}(fu)+(\tau
g)\pp_q(fu)$ (via $\pp_q(ab)=(\tau a)\pp_qb+(\pp_a)\tau^{-1}b$).
The first term in \eqref{4} becomes then
\bq\label{5}
-q\int u(\tau f)\pp_qg=-q\int \pp_qg\tau^{-1}(\tau u\tau^2 f)=
\end{equation}
$$=q\int \tau g\pp_q
(\tau u \tau^2f)=\int g\tau^{-1}\pp_q(\tau u\tau^2f)$$
Consequently \eqref{4} becomes
\bq\label{6}
\int g[\tau^{-1}\pp_q(\tau u\tau^2f)+\tau^{-1}(u\pp_qf)+aqc\pp_q^2(\tau+\tau^{-1})^{-1}
\pp_qf]
\end{equation}
Putting $f=u$ we obtain the Euler equation as in \cite{c1,c2}, namely
\bq\label{7}
qu_t=-qca\pp_q^2(\tau+\tau^{-1})^{-1}\pp_qu-\tau^{-1}\pp_q(\tau u\tau^2u)-
\tau^{-1}(u\pp_qu)
\end{equation}
Using also $\tau^{-1}\pp_q=q\pp_q\tau^{-1}$ we obtain \eqref{2} with $c'=ac$.  {\bf QED}
\end{theorem}
\indent
{\bf REMARK 3.1.}
One notes that \eqref{2} is morally equivalent to (4.19) in Theorem 4.3 of \cite{c1}
since $\psi''\sim \phi$ in the more tightly constrained framework of \cite{c1}.  Further
in view of the expression ${\bf (A14)}\,\,(\tau+\tau^{-1})^{-1}\sim \tau\sum
(-1)^n\tau^{2n}$ the equation \eqref{2} involves an infinite number of terms (much as
are indicated for qKdV in the hierarchy picture in \cite{c1,c2}.  Since we now have a
derivation with all of the classical algebraic and geometrical structure duplicated
it seems that \eqref{2} could be a good candidate for a canonical form.
Since the hierarchy qKdV is surely equally canonical one could anticipate an equivalence.
However the KdV equation arises in many different ways in mathematics and physics and
some caution should be employed in asserting that any particular form is canonical.
$\hfill\bs$


\end{document}